\documentclass[11pt]{amsart}

\usepackage[cp1250]{inputenc}
\usepackage{graphicx}
\usepackage{url}
\usepackage[all]{xy}
\usepackage{multicol}
\usepackage{amsthm}
\usepackage{amsmath}
\usepackage{amssymb}
\usepackage{geometry}
\usepackage{pinlabel}
\usepackage{hyperref}
\usepackage{array}
\usepackage{xcolor}
\usepackage{mathtools}

\newtheorem{Theorem}{Theorem}[section]
\newtheorem{lemma}[Theorem]{Lemma}
\newtheorem{proposition}[Theorem]{Proposition}
\newtheorem{corollary}[Theorem]{Corollary}
\newtheorem{definition}[Theorem]{Definition}

\newcounter{example}[section]
\newenvironment{example}[1][]{\refstepcounter{example}\par\medskip
   \noindent \textbf{Example~\theexample. #1} \rmfamily}{\medskip}

\newenvironment{remark}[1][]{\par\medskip
   \noindent \textbf{Remark.} \rmfamily}{\medskip}

\newcommand{\RR}{\mathbb {R}}
\newcommand{\CC}{\mathbb {C}}

\newcommand{\cl}{\mathop{\rm{cl}}}
\newcommand{\conv}{\mathop{\textrm{conv}\,}}
\newcommand{\ind}{\mathop{\rm{ind}}}
\newcommand{\inter}{\mathop{\rm{int}}}

\makeatletter
\@namedef{subjclassname@2020}{%
  \textup{2020} Mathematics Subject Classification}
\makeatother

\providecommand{\keywords}[1]
{
  \small	
  \textbf{Keywords :} #1 
}


\begin{document}

\title{Nonsmooth manifold decompositions}

\author{Eva Horvat}
\address{University of Ljubljana, Faculty of Education, Kardeljeva plo\v s\v cad 16, 1000 Ljubljana, Slovenia, eva.horvat@pef.uni-lj.si}

\begin{abstract}We study the structure induced on a smooth manifold by a continuous selection of smooth functions. In case such selection is suitably generic, it provides a stra\-ti\-fi\-cation of the manifold, whose strata are algebraically defined smooth submanifolds. When the continuous selection has nondegenerate critical points, stratification descends to the local topological structure. We analyze this structure for the maximum of three smooth functions on a 4-manifold, which provides a new perspective on the theory of trisections. 
\end{abstract}

\maketitle

\section{Introduction} \label{sec1}
Classical Morse theory is an indispensable tool in the study of low dimensional manifolds. By looking at the fibers of a smooth function, a manifold may be described as a union of fairly simple pieces. Assembling those pieces into few homeomorphic handlebo\-dies yields special decompositions: Heegaard splittings of 3-manifolds, trisections of 4-manifolds \cite{GK} and multisections of higher dimensional manifolds \cite{ACGK}.\\

A Heegaard decomposition of a closed, connected and oriented 3-manifold $Y$ is obtained from a classical Morse function $f\colon Y\to B^1$. Gay and Kirby have shown that a trisection of a closed, connected and oriented 4-manifold $X$ may similarly be constructed from a suitable Morse 2-function $f\colon X\to B^{2}$ \cite{GK}. On the other hand, they described how the handle decomposition, corresponding to a classical Morse function on a 4-manifold, fits into a trisection. Though natural and illuminating, neither of the two perspectives is quite successful in explaining the triple symmetry of a trisected manifold. Continuous selections of smooth functions offer an alternative between the 2-dimensional and the classical 1-dimensional Morse functions. Instead of describing a manifold in a linear fashion -- from bottom up -- they allow its construction ``from inside out'' in several directions.  \\

Originally, trisections of closed 4-manifolds were defined in the smooth category. However, some amount of research has lately been focused on trisections of PL manifolds that respect the underlying (singular) triangulations \cite{BHRT, CC, ST}. In another spirit, we study trisections induced by a continuous selection of smooth functions on a smooth manifold. This point of view offers several advantages: the fibers are easier to grasp, the spine and the trisection surface are algebraically defined, while handle decomposition carries a triple symmetry due to the underlying stratification.\\

A continuous selection of smooth functions $f_{1},\ldots ,f_{m}$ on a manifold $M$ is a continuous function that coincides with at least one of the functions $f_{i}$ at each point of $M$. Applying generalized derivation, introduced by Clarke \cite{CL}, the topology of $M$ may be studied by a suitably adapted version of Morse theory. We begin by analyzing the structure, revealed by a continous selection of smooth functions on a manifold. Under some nondegeneracy conditions, such structure defines a stratification of the manifold, whose strata are algebraically defined smooth submanifolds. We show examples of 3-dimensional Heegaard splittings and 4-dimensional trisections as instances of such stratifications. Moreover, the local topological structure around the critical point of a Morse CS function $f=\{f_{1},f_{2},f_{3}\}$ is studied, and a stratified handle decomposition, corresponding to such function, is des\-cribed.\\

The paper is organized as follows. In Section \ref{sec2}, we briefly present the basic concepts from nonsmooth analysis that we will need. In Section \ref{sec3} we study CS functions and show that the structure, induced on a smooth manifold by such function, defines its stratification. In Section \ref{sec4}, we focus on 4-manifolds. In Subsection \ref{subs41}, we recall the definition of trisection and describe CS functions that induce trisections of 4-manifolds. In Subsection \ref{subs42}, we investigate the local structure of a 4-manifold around a critical point of a CS Morse function $\max \{f_{1},f_{2},f_{3}\}$, which gives rise to stratified handles with triple symmetry.

\section*{Acknowledgements}

The author would like to thank Peter Feller for several discussions that led to the ideas presented in this paper. This research was supported by the Slovenian Research Agency grants P1-0292, J1-4031 and N1-0278.

\section{Preliminaries from nonsmooth analysis} \label{sec2}

A classical reference on nonsmooth analysis is Clarke's book \cite{CL}. Morse theory adapted to the context of piecewise smooth functions was studied in \cite{BO,JP,BC,BKS,APS}. Here we recall some basic definitions and results. For the remainder of this section, we denote by $M$ a smooth $n$-manifold. 

\begin{definition} Let $f\colon M\to \RR $ be a continuous function and let $f_{1},f_{2},\ldots ,f_{m}\colon M\to \RR $ be smooth functions. The function $f$ is called a \textbf{continuous selection} of the functions $f_{1},f_{2},\ldots ,f_{m}$ if it satisfies the following two conditions:
\begin{enumerate}
\item the set $I_{f}(x)=\left \{i\in \{1,2,\ldots ,m\}| f(x)=f_{i}(x)\right \}$ is nonempty for every point $x\in M$,
\item for any index $i\in \{1,2,\ldots ,m\}$, there exists a point $x\in M$ with $i\in I_{f}(x)$.
\end{enumerate}
Denote by $CS(f_{1},\ldots ,f_{m})$ the set of all continuous selections of the functions $f_{1}, f_{2}, \ldots ,f_{m}$. The set $I_{f}(x)$ is called the \textbf{active index set} of $f$ at the point $x$. Furthermore, the set $$\widehat{I}_{f}(x_0)=\left \{i\in \{1,\ldots ,m\}\, |\, x_{0}\in \cl (\inter \{x\, |\, f(x)=f_{i}(x)\})\right \}$$ is called the set of \textbf{essentially active indices} (where `` $\inter \,S$'' and `` $\cl S$'' denote the interior and the closure of a set $S$). A function $f_i$ is called \textbf{essentially active} at $x_0$ if $i\in \widehat{I}_{f}(x_{0})$. 
\end{definition}

It follows from the above definition that $\widehat{I}_{f}(x)\subseteq I_{f}(x)$ for every $x\in M$. For every piecewise differentiable function $f$ on $M$ and every $x\in M$ there exists a collection $\{f_{1},f_{2},\ldots ,f_{m}\}$ of smooth functions that are essentially active in a neighborhood of $x$ \cite[Proposition 4.1.1]{SCH}. Moreover, a continuous selection of smooth functions is locally Lipschitz continuous. In this setting, generalized differentiation defined by Clarke has the following form. 

\begin{definition} Let $f$ be a continuous selection of smooth functions $f_{1},\ldots ,f_{m}\colon \RR ^{n}\to \RR $. 
Clarke's  \textbf{subdifferential} of $f$ at a point $x \in \RR ^n$ is defined by $$\partial f(x)=\conv \left \{ \nabla f_{i}(x)\, |\, i\in \widehat{I}_{f}(x)\right \}\;,$$ where $\conv S$ denotes the convex hull of a set $S$. A point $x_{0}\in \RR ^n$ is called a \textbf{critical point} of $f$ if $0\in \partial f(x_{0})$.  
\end{definition}

More generally, when $f$ is a continuous selection of smooth functions on a smooth ma\-ni\-fold $M$, the Clarke's subdifferential of $f$ may be defined in a similar manner using the local coordinate charts on $M$. Let $\phi \colon U\to M$ be a local chart with $U\subset \RR ^{n}$ and $\phi (0)=x$, then we define $\partial f(x)=\conv \left \{ \nabla (f_{i}\circ \phi )(0)\, |\, i\in \widehat{I}_{f}(x)\right \}$. 

A version of the first Morse Lemma for locally Lipschitz continuous functions was proved in \cite{APS}.

\begin{proposition}[\cite{APS}] \label {prop0} Let $f\colon M\to \RR $ be a locally Lipschitz continuous function and denote by $M_{t}=f^{-1}((-\infty ,t])$ its lower level set for $t\in \RR $. Suppose $M_{\beta }$ is compact and $f^{-1}([\alpha ,\beta ])$ does not contain any critical point of $f$. Then there exists a Lipschitz continuous mapping $F\colon M_{\beta }\times [0,1]\to M_{\beta }$ such that $$F(M_{\beta },t)\subset M_{\beta +t(\alpha -\beta )}\,,\quad F(x,t)=x\; \textit{for all $x\in M_{\alpha }$ and all $t\in [0,1]$\;.}$$ 
\end{proposition}

In order to study the local behaviour of $M$ around a critical point of a piecewise diffe\-ren\-tiable function, we need the following definition. 

\begin{definition} \label{def1} Let $f_{1},\ldots ,f_{m}\colon M\to \RR $ be smooth functions and $f\in CS(f_{1},\ldots ,f_{m})$. A critical point $x_{0}\in M$ of $f$ is called \textbf{nondegenerate} if the following two conditions hold:
\begin{itemize}
\item[(ND1)] For each $i\in \widehat{I}_{f}(x_{0})$, the set of differentials $\left \{\nabla f_{j}(x_{0})\, |\, j\in \widehat{I}_{f}(x_{0})\backslash \{i\}\right \}$ is linearly independent;
\item[(ND2)] The second differential $\nabla _{x}^{2}L(x,\lambda )(x_0)$ is regular on $$\widehat{T}(x_0)=\bigcap _{i\in \widehat{I}_{f}(x_0)}Ker (\nabla f_{i}(x_0))\;,\quad \textrm{where}\quad L(x,\lambda )=\sum _{i\in \widehat{I}_{f}(x_0)}\lambda _{i}f_{i}(x)\;,$$
and the numbers $\lambda _{i}\in \RR $ are such that $$\sum _{i\in \widehat{I}_{f}(x_0)}\lambda _{i}\nabla f_{i}(x_{0})=0\,,\quad \sum _{i\in \widehat{I}_{f}(x_0)}\lambda _{i}=1\,,\quad \lambda _{i}\geq 0\textrm{ for every }i\in \widehat{I}_{f}(x_0)\;.$$
\end{itemize}
In this case, the \textbf{quadratic index} of the critical point $x_0$ is the dimension of a maximal linear subspace of $\widehat{T}(x_0)$ on which the quadratic form $y^{T}\nabla _{x}^{2}L(x,\lambda )(x_0)y$ is negative definite. 
\end{definition}

A version of the second Morse lemma for continuous selections of smooth functions was proved by Jongen and Pallaschke \cite{JP}. It requires the following definition. 

\begin{definition} Two continuous functions $f,g\colon M\to \RR $ are said to be \textbf{topologically equivalent} at $(x_{0},y_{0})\in M\times M$ if there exist open subsets $U$ of $x_0$ and $V$ of $y_0$ in $M$ and a homeomorphism $\phi \colon U\to V$, for which $\phi (x_0)=y_0$ and $f\circ \phi ^{-1}=g$ on $V$. 
\end{definition}

\begin{Theorem}[\cite{JP}] \label{th} Let $f_{1},\ldots ,f_{m}\colon M\to \RR $ be twice continuously differentiable functions on an $n$-manifold $M$ and let $f\in CS(f_{1},\ldots ,f_{m})$. Choose local coordinates $(y_{1},y_{2},\ldots ,y_{n})$ in the neighborhood of a point $x_{0}\in M$. Then the following holds: \begin{itemize}
\item[(i)] if $x_0$ is a noncritical point of $f$, then $f$ and $f(x_0)+y_1$ are locally topologically equivalent at the point $(x_{0},0)$;
\item[(ii)] if $x_0$ is a nondegenerate critical point of $f$, then $f$ is locally topologically equivalent at $(x_{0},0)$ to a function of the form $$f(x_0)+g(y_{1},\ldots ,y_{k})-\sum _{i=k+1}^{k+m}y_{i}^{2}+\sum _{j=k+m+1}^{n}y_{j}^{2}\;,$$ where $k=|\widehat{I}_{f}(x_0)|-1$, $g\in CS(y_{1},\ldots ,y_{k},-\sum _{i=1}^{k}y_{i})$, and $m$ is the quadratic index of $f$ at $x_0$.
\end{itemize}
\end{Theorem}

In particular, Theorem \ref{th} implies that the nondegenerate critical points of the function $f\in CS(f_{1},\ldots ,f_{m})$ are isolated.  

\begin{example} \label{ex4} Define functions $f_{1},f_{2}\colon \RR ^{3}\to \RR $ by $f_{1}(x,y,z)=x^{2}-y^{2}, f_{2}(x,y,z)=y^{2}-z^{2}$ and let $f$ denote the continuous selection $f=\max \{f_{1},f_{2}\}$. Consider the restriction of $f$ to the 2-sphere $S^{2}=\{(x,y,z)\in \RR ^{3}\,|\, x^{2}+y^{2}+z^{2}=1\}$. There are two disjoint circles of points $x\in S^{2}$ for which $\widehat{I}_{f}(x)=\{1,2\}$, namely $C_{\pm}=\left \{\left (x,\pm \frac{1}{\sqrt{3}},z\right )\in \RR ^{3}\, |\, x^{2}+z^{2}=\frac{2}{3}\right \}$. These two circles separate the sphere into an open annulus, where $I_{f}=\{1\}$, and two open disks on which $I_{f}=\{2\}$. 
\begin{figure}[h!]
\labellist
\normalsize \hair 2pt
\pinlabel $x$ at 57 42
\pinlabel $y$ at 170 73
\pinlabel $z$ at 90 182
\endlabellist
\begin{center}
\includegraphics[scale=0.63]{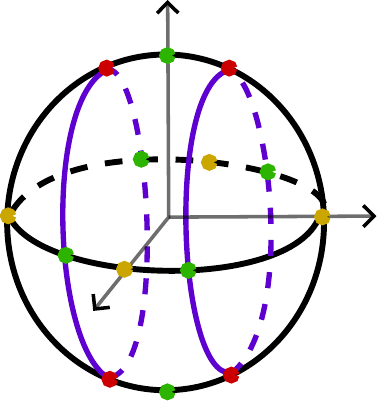}
\caption{Critical points of the function $f$ on $S^2$ from Example \ref{ex4} are indicated by red (quadratic index 0), green (quadratic index 1) and yellow (index 2) dots.}
\label{fig7}
\end{center}
\end{figure}
The annulus contains two critical points $(\pm 1,0,0)$ of index 2 and two critical points $(0,0,\pm 1)$ of index 1. Each of the disks with $I_{f}=\{2\}$ contains a critical point $(0,\pm 1,0)$ of index 2. To find the ``nonsmooth'' critical points of $f$ that lie on the circles $C_{+}$ and $C_{-}$, we express $f$ in the local coordinates $\phi _{\pm }\colon \RR ^{2}\to S^{2}$, given by $\phi _{\pm }(x,z)=(x,\pm \sqrt{1-x^{2}-z^{2}},z)$, and obtain $(f_{1}\circ \phi _{\pm })(x,z)=2x^{2}+z^{2}-1, (f_{2}\circ \phi _{\pm })(x,z)=1-x^{2}-2z^{2}$. Critical points are the solutions of 
\begin{xalignat*}{1}
t\,\nabla (f_{1}\circ \phi _{\pm })+(1-t)\nabla (f_{2}\circ \phi _{\pm })=(6t-2)x\frac{\partial }{\partial x}+(6t-4)z\frac{\partial }{\partial z}=0\;.
\end{xalignat*} At $t=\frac{2}{3}$, there are four nondegenerate critical points $\left (0,\pm \sqrt{\frac{1}{3}},\pm \sqrt{\frac{2}{3}}\right )$ of quadratic index 0. At each of these points, $f$ is locally topologically equivalent to the function $-\frac{1}{3}+|y|+x^{2}$. At $t=\frac{1}{3}$, there are four nondegenerate critical points $\left (\pm \sqrt{\frac{2}{3}},\pm \sqrt{\frac{1}{3}},0\right )$ of quadratic index 1. At each of these points, $f$ is locally topologically equivalent to the function $\frac{1}{3}+|y|-z^{2}$.
\end{example}



\section{Stratification of a manifold, induced by a CS function} \label{sec3}

A smooth function on a manifold reveals important information about its basic constituent parts. What happens if we observe a continuous selection of smooth functions instead? It turns out that ``corners'' of various dimensions provide an additional structure. \\

\subsection{CS functions on manifolds}

In this paper, we will focus on the following two types of functions. 

\begin{definition} \label{def2} Let $M$ be a smooth manifold. A \textbf{CS function} on $M$ is a function $f\colon M\to \RR $ for which there exist smooth functions $f_{1},\ldots ,f_{m}\colon M\to \RR $ such that \begin{enumerate}
\item $f\in CS(f_{1},f_{2},\ldots ,f_{m})$,
\item at any point $x\in M$ we have $I_{f}(x)=\widehat{I}_{f}(x)$,
\item at any point $x\in M$, the gradients of the active functions $\{f_{j}\,|\, j\in I_{f}(x)\}$ are affinely independent.
\end{enumerate} A \textbf{CS Morse function} on $M$ is a CS function on $M$ whose all critical points are nondegenerate. 
\end{definition}

Throughout this section, we denote by $M$ a smooth $n$-manifold and by $f$ a CS function on $M$, given as the continuous selection of smooth functions $f_{1},f_{2},\ldots ,f_{m}\colon M\to \RR $. For any subset $J\subset \{1,2,\ldots ,m\}$, we denote by $$M_{J}=\left \{x\in M\, |\, I_{f}(x)=J\right \}$$ the subset of points where the active index set equals $J$, and let $N_{J}=\{x\in M\, |\, J\subseteq I_{f}(x)\}$. Note that $M_{J}\subseteq N_{J}$ for any $J\subset \{1,2,\ldots ,m\}$. 

\begin{lemma} \label{lemma1} $M_{\{i\}}$ is an open submanifold of $M$ for any $i\in \{1,2,\ldots ,m\}$. Moreover, $N_{J}$ is a closed subset of $M$ for any $J\subseteq \{1,2,\ldots ,m\}$. 
\end{lemma}
\begin{proof} Let us denote $M_{f_i}=\left \{x\in M\, |\, f_{i}(x)=f(x) \right \}$. Since $f_{i}- f$ is continuous, the set $M_{f_i}=(f_{i}-f)^{-1}(0)$ is a closed subset of $M$. Its complement $M\backslash M_{f_i}$ is an open subset of $M$ for every $i\in \{1,2,\ldots ,m\}$, and $$M_{\{i\}}=\bigcap _{j\in \{1,2\ldots ,m\}\backslash \{i\}}(M\backslash M_{f_{j}})$$ is an open subset of $M$ (thus a submanifold). For any collection of indices $J\subseteq \{1,2,\ldots ,m\}$, the set $N_{J}=\left \{x\in M\, |\, J\subseteq I_{f}(x)\right \}=\bigcap _{j\in J}M_{f_{j}}$ is a closed subset of $M$.
\end{proof}

It follows from Lemma \ref{lemma1} that $M_{\{i\}}$ is a smooth submanifold of codimension 0 in $M$. More generally, each $M_J$ is a submanifold of $M$. To prove this, we need the following standard result.  

\begin{proposition}[\textup{\cite[Proposition 5.28]{LE}}] \label{prop1} Let $S$ be a subset of a smooth $n$-manifold $M$. Then $S$ is an embedded $k$-submanifold of $M$ if and only if every point $p\in S$ has a neighborhood $U$ in $M$ such that $U\cap S$ is a level set of a submersion $F\colon U\to \RR ^{n-k}$. 
\end{proposition}

\begin{proposition} \label{prop2} Let $J\subset \{1,2,\ldots ,m\}$ be any subset such that $I_{f}(x)=J$ for some $x\in M$. Then $M_J$ is a smooth submanifold of codimension $|J|-1$ in $M$.  
\end{proposition}
\begin{proof} For $|J|=1$, the statement follows from Lemma \ref{lemma1}. 

Suppose that $J=\{i,j\}$ for two indices $i\neq j\in \{1,2,\ldots ,m\}$ and that $M_J\neq \emptyset $. Choose an arbitrary $x\in M_J$. Let $U$ be a neighborhood of $x$ in $M$ such that for any index $r\in \{1,2,\ldots ,m\}\backslash J$, we have $U\cap N_{\{r\}}=\emptyset $ (such a neighborhood exists since $N_{\{r\}}$ is a closed subset of $M$ that does not contain $x$ by Lemma \ref{lemma1}). Define a map $F\colon U\to \RR $ by $F(w)=f_{i}(w)-f_{j}(w)$. Then $F$ is a smooth function on $U$ and since at any point $w\in U$ the gradients $\nabla f_{i}(w)$ and $\nabla f_{j}(w)$ are affinely independent, $F$ is a submersion. Since $U\cap M_{\{i,j\}}=F^{-1}(0)$, Proposition \ref{prop1} implies that $M_{\{i,j\}}=M_J$ is an embedded smooth submanifold of codimension 1 in $M$.    

Now let $J=\{i_{1},i_{2},\ldots ,i_{k}\}\subset \{1,2,\ldots ,m\}$ be any subset with $M_J\neq \emptyset $ and choose an arbitrary $x\in M_J$. Let $U$ be a neighborhood of $x$ in $M$ such that for any index $r\in \{1,2,\ldots ,m\}\backslash J$, we have $U\cap N_{\{r\}}=\emptyset $. Define a map $F\colon U\to \RR ^{k-1}$ by $$F(w)=\left (f_{i_{1}}(w)-f_{i_{k}}(w),f_{i_{2}}(w)-f_{i_{k}}(w),\ldots ,f_{i_{k-1}}(w)-f_{i_{k}}(w)\right )\;.$$ At any point $w\in U$, the set $\{\nabla f_{i_{j}}(w)\,|\, j\in \{1,2,\ldots ,k\}\}$ is affinely independent, thus $F$ is a submersion. Since $U\cap M_{J}=F^{-1}(0)$, Proposition \ref{prop1} implies that $M_J$ is an embedded smooth submanifold of codimension $k-1$ in $M$.   
\end{proof}

For any $J\subseteq \{1,2,\ldots ,m\}$, we denote by $\partial M_{J}$ the topological boundary (frontier) of the subset $M_{J}\subset M$. Note that $M_{J}\cap (\partial M_J)=\emptyset $, since $M_{J}$ is a submanifold of $M$. 

\begin{lemma}\label{lemma1.1} For any $i\in \{1,2,\ldots ,m\}$, we have $$M_{\{1,2,\ldots ,m\}}=\partial M_{\{1,2,\ldots ,m\}\backslash \{i\}}=\bigcap _{j=1}^{m}\partial M_{\{1,2,\ldots ,m\}\backslash \{j\}}\;.$$  
\end{lemma}
\begin{proof} Recall that $N_{\{1,2,\ldots ,m\}\backslash \{i\}}=\{x\in M\,|\, \{1,2,\ldots ,m\}\backslash \{i\}\subseteq I_{f}(x)\}$. By Lemma \ref{lemma1}, $N_{\{1,2,\ldots ,m\}\backslash \{i\}}$ is a closed subset of $M$ that contains $M_{\{1,2,\ldots ,m\}\backslash \{i\}}$. Moreover, $$M_{\{1,2,\ldots ,m\}\backslash \{i\}}=\{x\in N_{\{1,2,\ldots ,m\}\backslash \{i\}}\,|\, f(x)\neq f_{i}(x)\}=N_{\{1,2,\ldots ,m\}\backslash \{i\}}\cap (f-f_{i})^{-1}(\RR \backslash \{0\})$$ is an open subset of $N_{\{1,2,\ldots ,m\}\backslash \{i\}}$ whose closure equals $N_{\{1,2,\ldots ,m\}\backslash \{i\}}$, and its boundary is exactly $\{x\in N_{\{1,2,\ldots ,m\}\backslash \{i\}}\,|\, f(x)=f_{i}(x)\}=M_{\{1,2,\ldots ,m\}}$. It follows that $\partial M_{\{1,2,\ldots ,m\}\backslash \{i\}}=\partial M_{\{1,2,\ldots ,m\}\backslash \{j\}}$ for any $i,j\in \{1,2,\ldots ,m\}$. 
\end{proof}

\begin{lemma} \label{lemma2} For any $i\in \{1,2,\ldots ,m\}$, we have $M_{\{i\}}=\inter N_{\{i\}}$. 
\end{lemma}
\begin{proof} By Lemma \ref{lemma1}, $M_{\{i\}}$ is an open subset of $M$ that is clearly contained in $N_{\{i\}}$, thus $M_{\{i\}}\subseteq \inter N_{\{i\}}$. 

To show the other inclusion, first observe that $\{M_J\,|\, J\subseteq \{1,2,\ldots ,m\}\}$ is a partition of $M$ into pairwise disjoint smooth submanifolds by Proposition \ref{prop2}. Suppose that $\inter N_{\{i\}}\cap \inter N_{\{j\}}\neq \emptyset $ for some $j\neq i$, then this intersection is an open subset of $M$ that intersects at least one of the submanifolds $M_{\{r\}}$ of codimension 0 nontrivially (in fact, in an open set), thus $f_{i}(z)=f_{j}(z)=f(z)$ for some $z\in M_{\{r\}}$, which gives a contradiction. Thus, the subsets $\inter N_{\{r\}}$ are pairwise disjoint for $r\in \{1,2,\ldots ,m\}$. Now choose any $y\in \inter N_{\{i\}}$, and let $U$ be a neighborhood of $y$ in $M$ such that $U\subset \inter N_{\{i\}}$. Suppose that $f_{j}(y)=f(y)$ for some $j\in \{1,2,\ldots ,m\}\backslash \{i\}$. It follows from our reasoning above that $y\notin \inter N_{\{j\}}$, therefore $y\in N_{\{j\}}\backslash (\inter N_{\{j\}})$. By the Condition (2) of Definition \ref{def2}, we have $j\in I_{f}(y)=\widehat{I}_{f}(y)$ and thus $y\in \cl (\inter N_{\{j\}})$. But then the neighborhood $U$ intersects $\inter N_{\{j\}}$, which gives a contradiction. It follows that $I_{f}(y)=\{i\}$ and thus $y\in M_{\{i\}}$.  
\end{proof}

\begin{lemma} \label{lemma3} For any $J\subseteq \{1,2,\ldots ,m\}$ we have $$\partial M_{J}=\bigcup _{J\subsetneq I}M_{I}\;.$$
\end{lemma}
\begin{proof} Suppose that $x\in \partial M_J$. For every $j\in J$, we have $f_{j}|_{M_{J}}=f|_{M_{J}}$ and it follows by continuity that $f_{j}(x)=f(x)$. Thus, $J\subset I_{f}(x)$ and since $x\notin M_J$, $J$ is a proper subset of $I_{f}(x)$. This shows that $\partial M_{J}\subseteq \bigcup _{J\subsetneq I}M_{I}$. 

Let us show the other inclusion by induction on $|J|$. To begin, let $J=\{j\}$ and choose an element $x\in M$ with $J\subsetneq I_{f}(x)$. By Condition (2) of Definition \ref{def2}, we have $j\in \widehat{I}_{f}(x)$ and thus $x\in \cl (\inter N_{\{j\}})=\cl (M_{\{j\}})$ by Lemma \ref{lemma2}. Since $x\notin M_{\{j\}}$, it follows that $x\in \partial M_{\{j\}}$. 

Now suppose that for some $n\geq 1$ and every $J\subseteq \{1,2,\ldots ,m\}$ with $|J|=n$ we have $\bigcup _{J\subsetneq I}M_{I}\subseteq \partial M_J$. Choose any subset $J'\subseteq \{1,2,\ldots ,m\}$ with $|J'|=n+1$ and write $J'=J\sqcup \{j\}$. By Proposition \ref{prop2}, $M_{J}$ is a submanifold of codimension $n-1$ and $M_{\{j\}}$ is a submanifold of codimension 0. Thus, $\partial M_{J}$ is a closed topological submanifold of codimension $n$ and $\partial M_{\{j\}}$ is a closed topological submanifold of codimension 1. It follows by the induction hypothesis that $\partial M_{J}=\bigcup _{J\subsetneq I}M_{I}$ and $\partial M_{\{j\}}=\bigcup _{\{j\}\subsetneq I}M_{I}$, therefore 
\begin{xalignat}{1}\label{dec}
& \partial M_{J}\cap \partial M_{\{j\}}=\bigcup _{J'\subseteq I}M_{I}=M_{J'}\cup \bigcup _{J'\subsetneq I}M_{I}
\end{xalignat} is a closed topological submanifold of codimension $n$. In the decomposition \eqref{dec}, $M_J'$ is a submanifold of codimension $n$ and each $M_I$ with $J'\subsetneq I$ is a submanifold of codimension $\geq n+1$. Thus, for every $x\in \bigcup _{J'\subsetneq I}M_{I}$ and every neighborhood $U$ of $x$ in $M$, we have $U\cap M_{J'}\neq \emptyset $. It follows that $\bigcup _{J'\subsetneq I}M_{I}\subseteq \cl (M_{J'})$ and since $\left (\bigcup _{J'\subsetneq I}M_{I}\right )\cap M_{J'}=\emptyset$, we have $\bigcup _{J'\subsetneq I}M_{I}\subseteq \partial M_{J'}$. 
\end{proof}

\begin{proposition} \label{prop3} For any $x\in M$ and any $J\subseteq \{1,2,\ldots ,m\}$, the following statements are equivalent: \begin{itemize}
\item[(a)] $x\in \cl (M_J)$,
\item[(b)] $f(x)=f_{j}(x)$ for every $j\in J$,
\item[(c)] $J\subseteq \widehat{I}_{f}(x)$,
\item[(d)] for every neighborhood $U$ of $x$ and every $j\in J$ we have $U\cap M_{\{j\}}\neq \emptyset $.
\end{itemize}
\end{proposition}
\begin{proof}
(a) $\Rightarrow $ (b) Choose an $x\in \cl (M_{J})$. Since $f|_{M_{J}}=f_{j}|_{M_J}$ for every $j\in J$, it follows by continuity that $f(x)=f_{j}(x)$ for every $j\in J$. \\ 
(b) $\Rightarrow $ (a) Suppose that $f(x)=f_{j}(x)$ for every $j\in J$, then $J\subseteq I_{f}(x)$. If $I_{f}(x)=J$, it follows that $x\in M_{J}$. Otherwise, we have $J\subsetneq I_{f}(x)$ and then Lemma \ref{lemma3} implies that $x\in \partial M_{J}$. \\  
(b) $\Rightarrow $ (d) Suppose that $f(x)=f_{j}(x)$ for every $j\in J$ and let $U$ be a neighborhood of $x$ in $M$. Then $J\subseteq I_{f}(x)=\widehat{I}_{f}(x)$ and it follows by Lemma \ref{lemma2} that $x\in \cl (\inter (N_{\{j\}}))=\cl (M_{\{j\}})$, therefore $U\cap M_{\{j\}}\neq \emptyset $. \\
(d) $\Rightarrow $ (b) Suppose that every neighborhood of $x$ intersects $M_{\{j\}}$ for every $j\in J$. Then $x\in \cl (M_{\{j\}})$ and since $f|_{M_{\{j\}}}=f_{j}|_{M_{\{j\}}}$, it follows by continuity that $f(x)=f_{j}(x)$ for every $j\in J$.\\
(c) $\Rightarrow $ (d) Let $J\subseteq \widehat{I}_{f}(x)$ and choose a neighborhood $U$ of $x$ in $M$. Then $x\in \cl (\inter N_{\{j\}})$ for every $j\in J$ and it follows by Lemma \ref{lemma2} that $U\cap M_{\{j\}}\neq \emptyset $. \\
(d) $\Rightarrow $ (c) Suppose that every neighborhood of $x$ intersects $M_{\{j\}}$ for every $j\in J$; then $x\in \cl (M_{\{j\}})=\cl (\inter N_{\{j\}})$, therefore $j\in \widehat{I}_{f}(x)$ for every $j\in J$.
\end{proof}

\subsection{Stratification}

We have shown that a CS function $f\in CS(f_{1},\ldots ,f_{m})$ on a smooth manifold $M$ induces a decomposition of $M$ into a union $$M=\bigcup _{J\subseteq \{1,2,\ldots ,m\}} M_J$$ of pairwise disjoint smooth submanifolds. Moreover, the closures of these submanifolds carry an extra structure, described in Lemma \ref{lemma3} and Proposition \ref{prop3}. This is an instance of a well-known and studied phenomenon, called stratification. 

\begin{definition}\cite{FR} A \textbf{filtered space} is a Hausdorff topological space $X$, equipped by a sequence of closed subspaces $$\emptyset =X^{0}\subset X^{1}\subset X^{2}\subset \ldots \subset X^{n-1}\subset X^{n}=X$$ for some integer $n\geq 0$. The space $X^i$ is called the $i$-\textbf{skeleton}, and the index $i$ is called the \textbf{formal dimension of the skeleton}. The connected components of $X_{i}=X^{i}-X^{i-1}$ are called the \textbf{strata} of $X$ of formal dimension $i$. If $X$ is a filtered space of formal dimension $n$, then the components of $X_{n}=X^{n}-X^{n-1}$ are called the \textbf{regular strata} of $X$ and all the other strata are called \textbf{singular strata}. 
\end{definition}

\begin{definition}\cite{FR} Let $X$ be a filtered space. The filtered set $X$ satisfies the \textbf{frontier condition} if for any two strata $S,T$ of $X$ the following holds:
$$\textit{ if $S\cap \cl (T)\neq  \emptyset $, then $S\subset \cl (T)$\;.}$$ A \textbf{stratified space} is a filtered space that satisfies the frontier condition.
\end{definition}

Consider our manifold $M$, endowed by the CS function $f\in CS(f_{1},\ldots ,f_{m})$, and define the skeleta of $M$ by $$X^{i}=\left \{x\in M\,|\, |I_{f}(x)|\geq m-i+1\right \}$$ for $0\leq i\leq m$.
\begin{proposition} The sequence $$\emptyset = X^{0}\subset X^{1}\subset \ldots \subset X^{m-1}\subset X^{m}=M$$ defines a stratification of $M$.
\end{proposition}
\begin{proof} We have $X^{0}=\emptyset $ and $X^{m}=M$ by definition. More generally, it follows from Lemma \ref{lemma3} and Proposition \ref{prop3} that $X^{i}=\bigcup _{|J|=m-i+1} \cl (M_J)$ is a closed subspace of $M$ for every $i\in \{0,1,\ldots ,m\}$. The regular strata of $M$ are submanifolds $M_{\{i\}}$ for $i\in \{1,2,\ldots ,m\}$. More generally, we have $$X^{i}-X^{i-1}=\left \{x\in M\,|\, m-i+1\leq |I_{f}(x)|<m-i+2\right \}=\bigcup _{|J|=m-i+1}M_{J}\;.$$ Thus, any submanifold $M_{J}$ with $|J|=i$ represents a stratum of formal dimension $m-i+1$. Any two strata $M_{J}$ and $M_{J'}$ with $J\neq J'$ are disjoint, while Lemma \ref{lemma3} implies that the boundary of any stratum is a disjoint union of strata. If we have $S\cap \cl (T)\neq \emptyset $ for two strata $S$ and $T$, then it follows that either $S=T$ or $S\subset \partial T$. Therefore the above filtration of $M$ satisfies the frontier condition. 
\end{proof}


What can we say about the fibers of the CS function $f\colon M\to \RR $? Let us denote $Y^{i}=X^{i}\cap f^{-1}(t)$ for $i=0,1,\ldots ,m$. 

\begin{corollary} The sequence of subspaces $$\emptyset =Y^{0}\subset Y^{1}\subset \ldots \subset Y^{m-1}\subset Y^{m}=f^{-1}(t)$$ defines a stratification of the fiber $f^{-1}(t)$ for any $t\in f(M)$.
\end{corollary}
\begin{proof} Since $f^{-1}(t)$ is a closed subspace of $M$, also $Y^{i}$ is a closed subspace of $M$. We have $Y^{i}-Y^{i-1}=X_{i}\cap f^{-1}(t)=\bigcup _{|J|=m-i+1}\{x\in M_{J}\,|\, f(x)=t\}$, therefore the subsets $\{x\in M_{J}\, |\, f(x)=t\}$ represent the strata of formal dimension $m-|J|+1$. 

Let $S$ and $T$ be two strata with $S\cap \cl (T)\neq \emptyset $. By Lemma \ref{lemma3} and Proposition \ref{prop3}, there exist subsets $J,K\subset \{1,2,\ldots ,m\}$ with $S=\{x\in M\,|\, I_{f}(x)=J,f(x)=t\}$ and $\cl (T)=\{x\in M\,|\, K\subseteq I_{f}(x), f(x)=t\}$, therefore $$S\cap \cl (T)=\{x\in M\,|\, K\subseteq I_{f}(x)=J,f(x)=t\}\neq \emptyset \;.$$ It follows that either $J=K$ or $K\subsetneq J$, thus $S=T$ or $S\subseteq \partial T$. We have shown that the above filtration of $f^{-1}(t)$ satisfies the frontier condition. 
\end{proof}  

\begin{remark} \label{rem1} Observe that the restriction of a CS function $f$ to every stratum is a smooth function, since $f|_{M_{J}}=f_{j}|_{M_{J}}$ for any $j\in J$. We might thus think of a CS function as a ``piecewise smooth'' function, where the smoothness condition holds on every stratum of our stratification (instead of a simplex in a triangulation). 
\end{remark}

\begin{example} \label{ex0} Define a CS function on the 3-sphere $$S^{3}=\{(x_{1},x_{2},x_{3},x_{4})\in \RR ^{4}\,|\, x_{1}^{2}+x_{2}^{2}+x_{3}^{2}+x_{4}^{2}=1\}$$ by $f(x_{1},x_{2},x_{3},x_{4})=\max \{x_{1}^{2}+x_{2}^{2},x_{3}^{2}+x_{4}^{2}\}$. Its image $f(S^{3})=\left [\frac{1}{2},1\right ]$ contains two critical values; the critical fiber $f^{-1}\left (\frac{1}{2}\right )$ is homeomorphic to the 2-torus, while the critical fiber $f^{-1}(1)$ looks like a Hopf link. Every regular fiber 
\begin{xalignat*}{1}
& f^{-1}(t)=\left \{(x_{1},x_{2},x_{3},x_{4})\,|\, x_{1}^{2}+x_{2}^{2}=t, x_{3}^{2}+x_{4}^{2}=1-t\right \}\cup \\
& \left \{(x_{1},x_{2},x_{3},x_{4})\,|\, x_{3}^{2}+x_{4}^{2}=t, x_{1}^{2}+x_{2}^{2}=1-t\right \}
\end{xalignat*}
for $\frac{1}{2}<t<1$ is a link of two tori that represent the boundary of a regular neighborhood of the Hopf link. The regular strata of the stratification induced by $f$ are the interiors of two solid tori $M_{\{1\}}=\left \{(x_{1},x_{2},x_{3},x_{4})\in S^{3}\,|\, x_{1}^{2}+x_{2}^{2}>x_{3}^{2}+x_{4}^{2}\right \}$ and $M_{\{2\}}=\left \{(x_{1},x_{2},x_{3},x_{4})\in S^{3}\,|\, x_{1}^{2}+x_{2}^{2}<x_{3}^{2}+x_{4}^{2}\right \}$, whose common boundary is the singular stratum $M_{\{1,2\}}=f^{-1}\left (\frac{1}{2}\right )$. Thus, stratification induced by $f$ gives the Heegaard decomposition of the 3-sphere of genus 1. 
\end{example} 

\begin{example} Note that the function from Example \ref{ex0} is not a CS Morse function. On the other hand, the linear function $f(x_{1},x_{2},x_{3},x_{4})=\max \{x_{1},x_{2}\}$ on $S^{3}$ is a CS Morse function with 4 critical points $\left (-\frac{1}{\sqrt{2}},-\frac{1}{\sqrt{2}},0,0\right )$, $\left (\frac{1}{\sqrt{2}},\frac{1}{\sqrt{2}},0,0\right )$, $(1,0,0,0)$ and $(0,1,0,0)$. A regular fiber $f^{-1}(t)$ for $-\frac{1}{\sqrt{2}}<t<\frac{1}{\sqrt{2}}$ is a union of two disks, glued along their boun\-da\-ry, and the critical fiber $f^{-1}(\frac{1}{\sqrt{2}})$ is a wedge of two spheres. A regular fiber $f^{-1}(t)$ for $\frac{1}{\sqrt{2}}<t<1$ is a disjoint union of two spheres. The stratification, induced by $f$, gives a Heegaard decomposition of genus 0. 
\end{example}


\section{CS functions inducing trisections} \label{sec4}

\subsection{Trisections} \label{subs41}

Our main motivation for studying CS functions is to understand their role in the theory of trisections. Recall the following definition from \cite{GK}. 

\begin{definition} Let $M$ be a closed, connected, oriented 4-manifold. Given two integers $0\leq k\leq g$, a \textbf{$(g,k)$-trisection} of $M$ is a decomposition $M=X_{1}\cup X_{2}\cup X_{3}$ such that: \begin{itemize}
\item $X_{i}$ is a 4-dimensional 1-handlebody obtained by attaching $k$ 1-handles to one 0-handle for each $i\in \{1,2,3\}$, 
\item $X_{i}\cap X_{j}$ is a 3-dimensional genus $g$ handlebody for each $i\neq j\in \{1,2,3\}$ and
\item $X_{1}\cap X_{2}\cap X_{3}$ is a closed surface of genus $g$. 
\end{itemize}
The surface $X_{1}\cap X_{2}\cap X_{3}$ is called the \textbf{trisection surface}, while the union of the 3-dimensional handlebodies $\cup _{i\neq j\in \{1,2,3\}}(X_{i}\cap X_{j})$ is called the \textbf{spine} of the trisection. 
\end{definition}

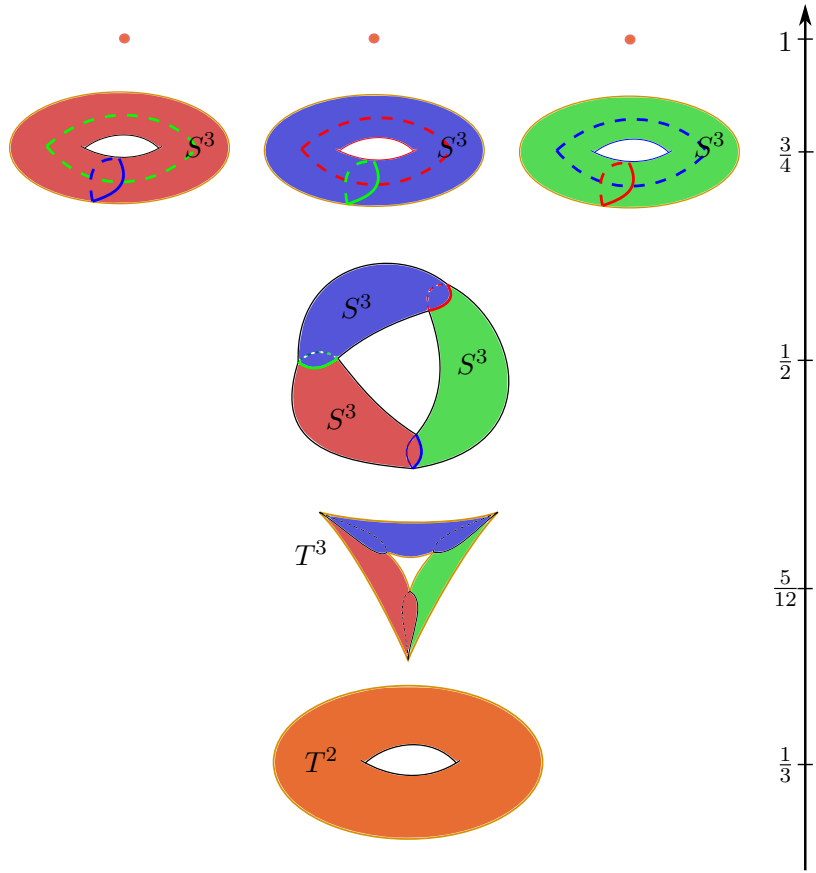
\begin{figure}[h!]
\labellist
\normalsize \hair 2pt
\pinlabel $\frac{1}{3}$ at 1540 210
\pinlabel $T^{2}$ at 620 220
\pinlabel $\frac{5}{12}$ at 1540 560
\pinlabel $\frac{1}{2}$ at 1540 1010
\pinlabel $T^3$ at 600 630
\pinlabel $\frac{3}{4}$ at 1540 1420
\pinlabel $S^3$ at 380 1440
\pinlabel $S^3$ at 880 1440
\pinlabel $S^3$ at 1390 1440
\pinlabel $1$ at 1540 1645
\pinlabel $S^3$ at 690 1120
\pinlabel $S^3$ at 660 900
\pinlabel $S^3$ at 920 1010
\endlabellist
\begin{center}
\includegraphics[scale=0.20]{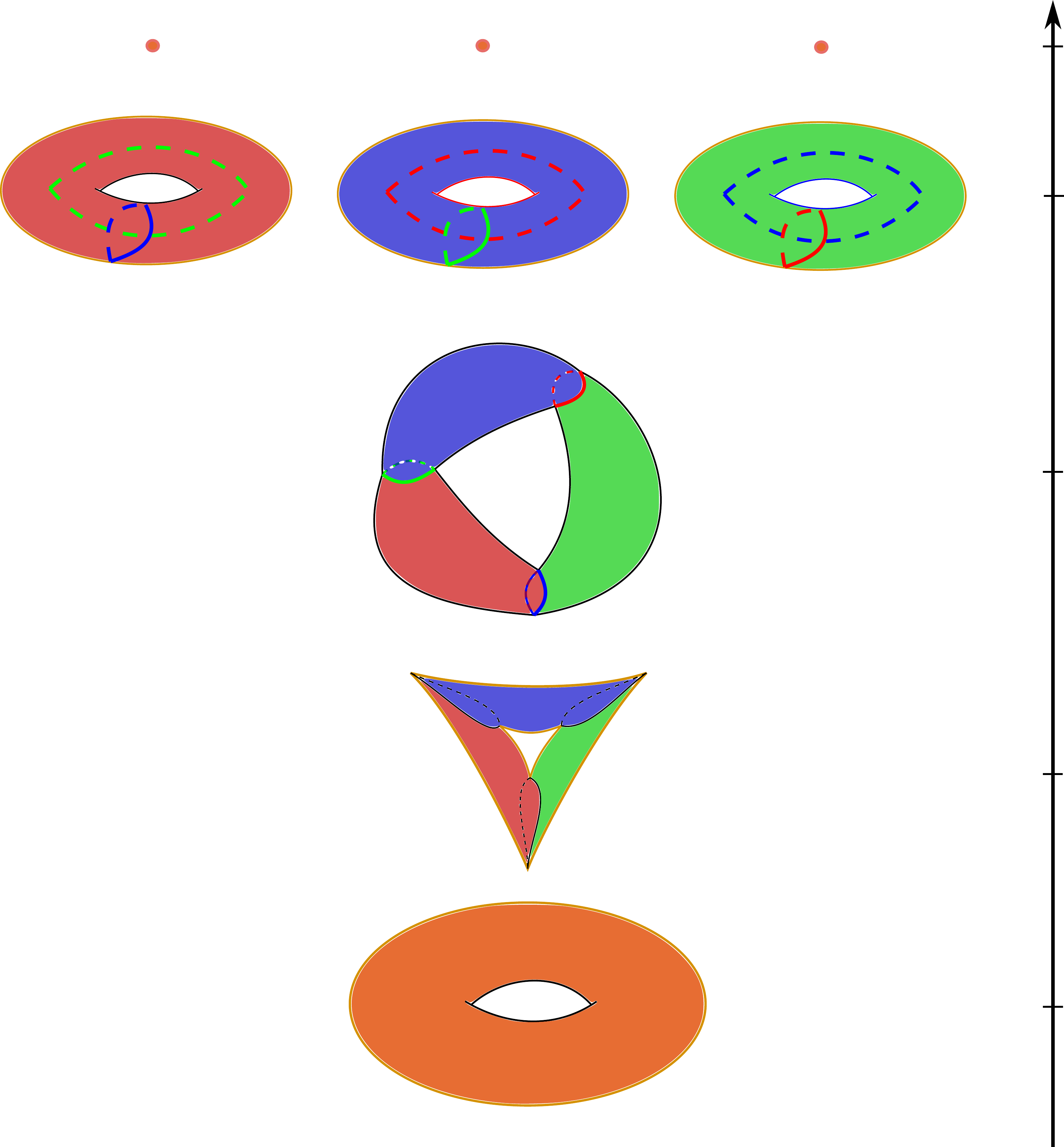}
\caption{The fibers of the CS function $f$ on $\CC P^2$ from Example \ref{ex3}.}
\label{fig2}
\end{center}
\end{figure}

\begin{example} \label{ex3} Define a CS function on the complex projective plane $\CC P^2$ by $$f([z_{0}:z_{1}:z_{2}])=\frac{1}{|z_{0}|+|z_{1}|+|z_{2}|}\max \left \{|z_{0}|,|z_{1}|,|z_{2}|\right \}\;.$$  The regular strata induced by $f$ are given by $$M_{\{i\}}=\{[z_{0}:z_{1}:z_{2}]\in \CC P^{2}\,|\, z_{i}=1,\, |z_{j}|< 1,\, |z_{k}|<1\textrm{ for $\{i,j,k\}=\{0,1,2\}$}\}\;,$$ each of which is diffeomorphic to $int(B^{2})\times int(B^{2})$. A singular stratum $$M_{\{i,j\}}=\{[z_{0}:z_{1}:z_{2}]\in \CC P^{2}\,|\, z_{i}=1,\, |z_{j}|=1,\, |z_{k}|<1\textrm{ for $\{i,j,k\}=\{0,1,2\}$}\}$$ is diffeomorphic to $S^{1}\times int(B^{2})$, and the final stratum $$M_{\{0,1,2\}}=\{[z_{0}:z_{1}:z_{2}]\in \CC P^{2}\,|\, z_{0}=1,\, |z_{1}|=|z_{2}|=1\}$$ is diffeomorphic to the torus. Thus, $f$ induces the standard (1,0)-trisection of $\CC P^{2}$.

The fibers of $f$ are sketched in Figure \ref{fig2}. The image $f(\CC P^{2})=\left [\frac{1}{3},1\right ]$ contains three critical values: $\frac{1}{3}$, $\frac{1}{2}$ and $1$. The restriction of $f$ to the regular strata contains one critical point in each stratum; these are the points $$[1:0:0]\in M_{\{0\}},\quad [0:1:0]\in M_{\{1\}},\quad [0:0:1]\in M_{\{2\}}\;,$$ whose union represents the critical fiber $f^{-1}(1)$. These three critical points are degenerate (the Hessian of $f_i$ is a singular matrix). Therefore, $f$ is not a CS Morse function. 

A regular fiber $f^{-1}(t)$ for $\frac{1}{2}<t<1$ is a disjoint union of three copies of $S^{3}$, given by 
$$\left \{[1\colon z\colon w]\, |\, |z|+|w|=\frac{1-t}{t}\right \}\cup \left \{[z\colon 1\colon w]\, |\, |z|+|w|=\frac{1-t}{t}\right \}\cup\left \{[z\colon w\colon 1]\, |\, |z|+|w|=\frac{1-t}{t}\right \}\;.$$ The first of these represents a 3-sphere with a Heegaard decomposition of genus 1 and Heegaard surface $\left \{[1\colon z\colon w]\, |\, |z|=|w|=\frac{1-t}{2t}\right \}$, while the circles $\left \{[1\colon z\colon 0]\, |\, |z|=\frac{1-t}{t}\right \}$ and $\left \{[1\colon 0\colon w]\, |\, |w|=\frac{1-t}{t}\right \}$ represent the cores of the two handlebodies. 

 The critical fiber $f^{-1}\left (\frac{1}{2}\right )$ is a union of three copies of $S^{3}$, glued together pairwisely along three circles of critical points. Each of the critical circles lives in a singular stratum $M_{\{i,j\}}$: $$\left \{[1\colon z\colon 0]\, |\, |z|=1\right \}\subset M_{\{0,1\}},\quad \left \{[z\colon 0\colon 1]\, |\, |z|=1\right \}\subset M_{\{0,2\}},\quad \left \{[0\colon 1\colon z]\, |\, |z|=1\right \}\subset M_{\{1,2\}}\;.$$ Each critical circle is the core of one of the 3-dimensional handlebodies that build the spine of the trisection. 
 
 A regular fiber $f^{-1}(s)$ for $\frac{1}{3}<s<\frac{1}{2}$ is a 3-torus $T^{3}$, given by \begin{xalignat*}{1}
 & \left \{[1\colon z\colon w]\, |\, |z|+|w|=\frac{1-s}{s}, |z|\leq 1, |w|\leq 1\right \}\cup \left \{[z\colon 1\colon w]\, |\, |z|+|w|=\frac{1-s}{s}, |z|\leq 1, |w|\leq 1\right \}\\
 & \cup  \left \{[z\colon w\colon 1]\, |\, |z|+|w|=\frac{1-s}{s}, |z|\leq 1, |w|\leq 1\right \}\;,\end{xalignat*} where each of the above subsets represents a product of a 2-torus with an interval, and they are glued together pairwisely along their boundary components. The critical fiber $f^{-1}\left (\frac{1}{3}\right )$ is a 2-torus, representing the singular stratum $M_{\{0,1,2\}}$, all of whose points are critical. 
\end{example}

Based on the above example, we formulate the following definition. 

\begin{definition} Let $M$ be a closed connected oriented 4-manifold and let $f\in CS(f_{1},f_{2},f_{3})$ be a CS function on $M$. We say that $f$ \textbf{induces a $(g,k)$ trisection} if $$\left (\cl (M_{\{1\}}),\cl (M_{\{2\}}),\cl (M_{\{3\}})\right )$$ defines a $(g,k)$ trisection of $M$. 
\end{definition}

The topology of a manifold $M$, revealed by a CS Morse function $f\in CS\{f_{1},\ldots ,f_{m}\}$, can also be seen through restrictions of $f$ to the submanifolds $M_J$ for $J\subset \{1,\ldots ,m\}$.   

\begin{lemma} \label{lemma4} Let $f\in CS(f_{1},f_{2},\ldots ,f_{m})$ be a CS Morse function. Then $f|_{M_{\{i\}}}=f_{i}|_{M_{\{i\}}}$ is a classical Morse function for every $i\in \{1,2,\ldots ,m\}$. 
\end{lemma}
\begin{proof} $M_{\{i\}}$ is an open submanifold of $M$ by Lemma \ref{lemma1} and the restriction $f|_{M_{\{i\}}}$ is smooth. For any critical point $x_0\in M_{\{i\}}$ we have $|\widehat{I}_{f}(x_0)|=1$, thus the nondegeneracy conditions (ND1) and (ND2) for $f$ from Definition \ref{def1} boil down to the classical nondegeneracy condition for $f_{i}|_{M_{\{i\}}}$. 
\end{proof}


\begin{lemma} \label{lemma5} Let $f\in CS(f_{1},f_{2},\ldots ,f_{m})$ be a CS Morse function on a manifold $M$. If $x_{0}\in M$ is a critical point of $f$ with $I_{f}(x_{0})=J$, then $x_{0}$ is a nondegenerate critical point of $f|_{M_J}$. 
\end{lemma}
\begin{proof} Suppose $x_0$ is a critical point of $f$ and $I_{f}(x_0)=J=\{f_{i_{1}},f_{i_{2}},\ldots ,f_{i_k}\}$. Let $\lambda _{1},\ldots ,\lambda _{k}\in [0,1]$ be such that $\sum _{j=1}^{k}\lambda _{j}=1$ and $\sum _{j=1}^{k}\lambda _{j}\nabla f_{i_{j}}(x_0)=0$. By Proposition \ref{prop2}, $M_{J}$ is a submanifold of $M$. Denote by $\iota \colon M_{J}\to M$ the associated embedding. Then 
$$\nabla \left (f_{i_{j}}|_{M_J}\right )(x_0)=\nabla (f_{i_{j}}\circ \iota )(x_0)=\textrm{proj}_{T_{x_{0}}M_{J}}\left (\nabla f_{i_{j}}(x_{0})\right )$$ for $j=1,2,\ldots ,k$. Here we denoted by $\textrm{proj}_{T_{x_{0}}M_{J}}(v)$ the projection of a tangent vector $v\in T_{x_{0}}M$ to the subspace $T_{x_{0}}M_{J}\leq T_{x_{0}}M$. It follows that $\sum _{j=1}^{k}\lambda _{j}\nabla \left (f_{i_{j}}|_{M_J}\right )(x_0)=0$, thus $x_0$ is a critical point of the restriction $f|_{M_J}$. Since $x_0$ is a nondegenerate critical point of $f$, it is a nondegenerate critical point of the restriction $f|_{M_J}$.
\end{proof}

The following theorem describes the general setting in which a CS function does induce a trisection of a 4-manifold. 

\begin{Theorem} \label{th1} Let $f\in CS(f_{1},f_{2},f_{3})$ be a CS Morse function on a closed, connected and oriented manifold $M$. Suppose that in each regular stratum $M_{\{i\}}$, $f$ has a single critical point of index 4, $k$ critical points of index 3 and no other critical points for $i\in \{1,2,3\}$. Moreover, suppose that in each singular stratum $M_{\{i,j\}}$, $f$ has a single critical point of index 3, $g$ critical points of index 2 and no other critical points for all $i\neq j\in \{1,2,3\}$. Then $f$ induces a $(g,k)$-trisection of the 4-manifold $M$. 
\end{Theorem}
\begin{proof} By Proposition \ref{prop2} and Lemma \ref{lemma1.1}, $M_{\{i,j\}}$ is a 3-manifold with boundary $M_{\{1,2,3\}}$. Observe that the restriction of $f$ to every stratum is a smooth function. Since $f$ is a CS Morse function, the restriction of $f$ to every stratum is a Morse function by Lemma \ref{lemma5}. 

The function $-f|_{M_{\{i,j\}}}$ has a single critical point of index 0 and $g$ critical points of index 1, thus $\cl (M_{\{i,j\}})=\cl (M_{\{i\}})\cap \cl (M_{\{j\}})$ is a 3-dimensional handlebody of genus $g$ and $M_{\{1,2,3\}}=\cl (M_{\{1\}})\cap \cl (M_{\{2\}})\cap \cl (M_{\{3\}})$ is a closed orientable surface of genus $g$. 

The function $-f|_{M_{\{i\}}}$ has a single critical point of index 0 and $k$ critical points of index 1, thus $\cl (M_{\{i\}})$ is a 4-dimensional handlebody of genus $k$ for $i=1,2,3$. 
\end{proof}

\subsection{The local structure of trisections, induced by a CS Morse function $\max \{f_{1},f_{2},f_{3}\}$} \label{subs42}

In this section, we investigate the local structure of a manifold around a nondegenerate critical point of a CS function. After some general observations, we focus on the case of a CS Morse function, inducing a trisection of a 4-manifold.  

The local topology of a manifold around a critical point in the case of the maximum (resp. minimum) function is particularly simple.

\begin{corollary}[\textup{\cite[Corollary 3.4]{JP}}] \label{cor3} Let $x_0$ be a nondegenerate critical point with quadratic index $r$ of $f\in CS (f_{1},f_{2},\ldots ,f_{m})$. Assume $f(x_0)=0$ and let $k=|\widehat{I}_{f}(x_{0})|-1$. \begin{itemize}
\item[(1)] If $f=\max \{f_{1},f_{2},\ldots ,f_{m}\}$, then $f$ is topologically equivalent to the function 
$$g(y_{1},y_{2},\ldots ,y_{n})=-\sum _{i=1}^{r}y_{i}^{2}+\sum _{j=r+1}^{n}y_{j}^{2}$$  at $(x_{0},0)$.
\item[(2)]  If $f=\min \{f_{1},f_{2},\ldots ,f_{m}\}$, then $f$ is topologically equivalent to the function $$g(y_{1},y_{2},\ldots ,y_{n})=-\sum _{i=1}^{r+k}y_{i}^{2}+\sum _{j=r+k+1}^{n}y_{j}^{2}$$ at $(x_{0},0)$. 
\end{itemize}
\end{corollary}

In the case of the maximum function, the index of a critical point is completely defined by the restriction of $f$ to the submanifold $M_J$. Denote by $\ind _{f}(x_0)$ the index of a critical point $x_0$ of the function $f$. 
 
\begin{proposition} \label{prop5} Let $f=\max \{f_{1},f_{2},\ldots ,f_{m}\}$ be a CS Morse function on a manifold $M$. If $x_0$ is a critical point of $f$ with $J=I_{f}(x_0)$, then $\ind_{f}(x_0)=\ind_{f|_{M_J}}(x_0)$.
\end{proposition}
\begin{proof} For $|J|=1$, the statement is obvious. 

Suppose that $J=\{i,j\}$ for some $i\neq j\in \{1,2,\ldots m\}$. By Proposition \ref{prop2}, $M_{\{i,j\}}$ is a codimension 1 submanifold of $M$. Since $x_{0}$ is a critical point of $f$, we have $0\in conv\{\nabla f_{i}(x_{0}),\nabla f_{j}(x_{0})\}$ and thus $t\nabla f_{i}(x_{0})+(1-t)\nabla f_{j}(x_0)=0$ for some $0\leq t\leq 1$. Since the gradient vectors $\nabla f_{i}(x_{0}),\nabla f_{j}(x_{0})$ are affinely independent, they span a 1-dimensional vector subspace of $T_{x_{0}}M$. When moving from the point $x_0$ in the direction of the gradient vector $\nabla f_{i}(x_{0})$, the value of $f_{i}$ increases and the value of $f_{j}$ decreases, so we end up at a point $x$ with $f_{i}(x)\neq f_{j}(x)$, thus $x\notin M_{\{i,j\}}$ (when moving from $x_0$ in the direction of $\nabla f_{j}(x_0)$, the conclusion is analogous). Therefore, the subspaces $Lin \{\nabla f_{i}(x_{0}),\nabla f_{j}(x_{0})\}$ and $T_{x_{0}}M_{\{i,j\}}$ have trivial intersection and consequently $T_{x_{0}}M=T_{x_{0}}M_{\{i,j\}}\oplus Lin \{\nabla f_{i}(x_{0}),\nabla f_{j}(x_{0})\}$. 

Choose a neighborhood $U_{x_0}$ of the critical point $x_0$ such that $U_{x_0}\cap \cl (M_{\{k\}})=\emptyset $ for $\{i,j,k\}=\{1,2,3\}$. It follows that $f|_{U_{x_0}}=\max \{f_{i},f_{j}\}=\frac{f_{i}+f_{j}}{2}+\frac{|f_{i}-f_{j}|}{2}$ and $f|_{U_{x_0}\cap M_{\{i,j\}}}=\frac{f_{i}+f_{j}}{2}$. The value of $f$ increases as we move away from $U_{x_0}\cap M_{\{i,j\}}$ in either direction of $\nabla f_{i}(x_0)$ or the direction of $\nabla f_{j}(x_0)$ . It follows from the argument above that all the negative eigenvalues of the second differential come from the restriction $f|_{M_{\{i,j\}}}$. 



Finally, let $J\subseteq \{1,2,\ldots ,m\}$ be any subset. Since $x_0$ is a critical point of $f$, we have $\sum _{j\in J}t_{j}\nabla f_{j}(x_0)=0$ for some $0\leq t_{j}\leq 1$ with $\sum _{j\in J}t_{j}=1$. By Proposition \ref{prop2}, $M_J$ is a submanifold of codimension $|J|-1$. By a similar argument as in the previous case, $Lin \{\nabla f_{j}(x_0)\,|\, j\in J\}\leq T_{x_{0}}M$ is a $(|J|-1)$-dimensional subspace whose intersection with $T_{x_{0}}M_{J}$ is trivial. Therefore, we obtain a decomposition $$T_{x_{0}}M=T_{x_{0}}M_{J}\oplus Lin \{\nabla f_{j}(x_{0})\,|\, j\in J\}$$ and any eigenvector of the second differential of $f$ corresponding to a negative eigenvalue lives in $T_{x_0}M_{J}\leq T_{x_{0}}M$. 
\end{proof}

For the remainder of this section, let $M$ be a closed connected 4-manifold and let $f=\max \{f_{1},f_{2},f_{3}\}$ be a CS Morse function on $M$. Suppose that $x_0$ is a critical point of $f$ with quadratic index $r$. Then Corollary \ref{cor3} implies there exists a neighborhood $U$ of $x_0$ in $M$, a neighborhood $V$ of $0\in \RR ^4$ and a homeomorphism $\phi \colon U\to V$ such that $\phi (x_0)=0$ and $f\circ \phi ^{-1}=g$, where $g(y_{1},y_{2},y_{3},y_{4})=-\sum _{i=1}^{r}y_{i}^{2}+\sum _{j=r+1}^{4}y_{j}^{2}$. Thus, the change in the topological type of the lower level sets at a critical point corresponds to the addition of a topological handle, and $f$ induces a topological handle decomposition of $M$. Stratification corresponding to the CS function $f$, however, endows this handle decomposition with an additional structure. 

\begin{definition} A \textbf{stratified $k$-handle} inside the stratified manifold $M$ is a topological handle $j\colon B^{k}\times B^{m-k}\to h\subset M$ such that
\begin{enumerate}
\item the core $j(B^{k}\times \{0\})$ is a $k$-dimensional disk, contained in a single stratum $M_J$, 
\item $h=\bigcup _{K\subseteq J}(h\cap M_{K})$,
\item the restriction $j|_{j^{-1}(h\cap M_{K})}$ is smooth for every $K\subseteq J$.    
\end{enumerate} A \textbf{stratified handle decomposition} of the stratified manifold $M$ is a sequence of subsets $$\emptyset =M^{(-1)}\subset M^{(0)}\subset M^{(1)}\subset M^{(2)}\subset M^{(3)}\subset M^{(4)}=M\;,$$ where each $M^{(k)}$ is obtained from $M^{(k-1)}$ by the attachement of stratified $k$-handles. 
\end{definition}

Let us show that the CS Morse function $f=\max \{f_{1},f_{2},f_{3}\}$ induces a stratified handle decomposition of the 4-manifold $M$. For each regular value of $f$, the level set $\{x\in M\, |\, f(x)=t\}$ is a closed 3-manifold, trisected into $\{x\in M\, |\, f_{i}(x)=t\}$ for $i=1,2,3$. Now let $x_{0}\in M$ be a critical point of $f$ with a 4-ball neighborhood $U_{x_0}$. At $x_0$, the topology of the lower level set $\{x\in M\, |\, f(x)\leq t\}$ changes by the addition of a topological handle. By Lemma \ref{lemma5} and Proposition \ref{prop5}, $x_0$ is also a critical point of $f|_{M_J}$ and $\ind _{f}(x_0)=\ind _{f|_{M_J}}(x_0)$, where $J=I_{f}(x_0)$. In other words, the core of the topological handle is a subset of the stratum $M_J$. By Lemma \ref{lemma3}, we have $M_{J}\subset \cl (M_{K})\Leftrightarrow K\subseteq J$, thus $U_{x_0}\cap M_{K}\neq \emptyset \Leftrightarrow K\subseteq J$. Since $f|_{M_K}=f_{i}|_{M_K}$ for any $i\in K$ and $f_{i}|_{M_K}$ is smooth, the restriction of the topological handle to each stratum $M_{K}$ with $K\subseteq J$ is smooth. Thus we obtain a stratified handle, that we will call the \textbf{stratified handle, induced by the} CS \textbf{Morse function} $f$. 

\begin{figure}[h!]
\labellist
\normalsize \hair 2pt
\pinlabel $x_0$ at 320 250
\pinlabel $M_{\{i,j\}}$ at 100 300
\pinlabel $M_{\{i,j\}}$ at 550 300
\pinlabel $M_{\{i\}}$ at 320 100
\pinlabel $M_{\{j\}}$ at 320 500
\pinlabel $x_0$ at 1080 250
\pinlabel $D$ at 1610 300
\pinlabel $x_0$ at 1900 290
\pinlabel $D$ at 2040 300
\pinlabel $D$ at 1140 310
\pinlabel $M_{\{1\}}$ at 1000 440
\pinlabel $M_{\{3\}}$ at 1280 310
\pinlabel $M_{\{2\}}$ at 1030 150
\endlabellist
\begin{center}
\includegraphics[scale=0.18]{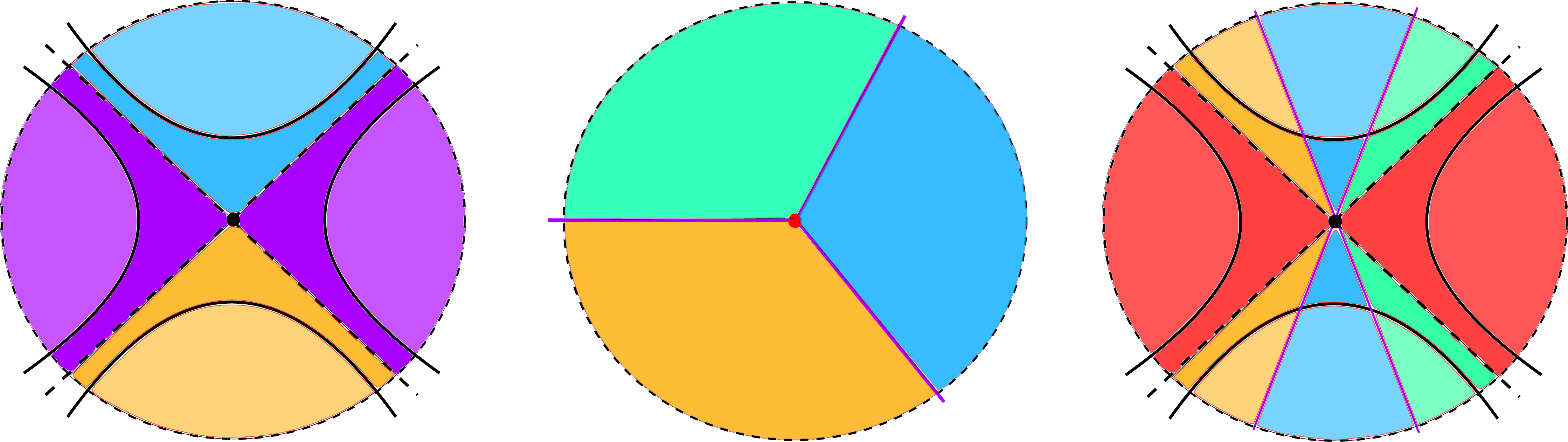}
\caption{The local structure of $M$ around a critical point $x_0$: \quad \quad \quad \quad (left) $x_{0}\in M_{\{i,j\}}$, (middle) $x_{0}\in M_{\{1,2,3\}}$, case (a), (right) $x_{0}\in M_{\{1,2,3\}}$, cases (b) and (c).}
\label{fig5}
\end{center}
\end{figure}
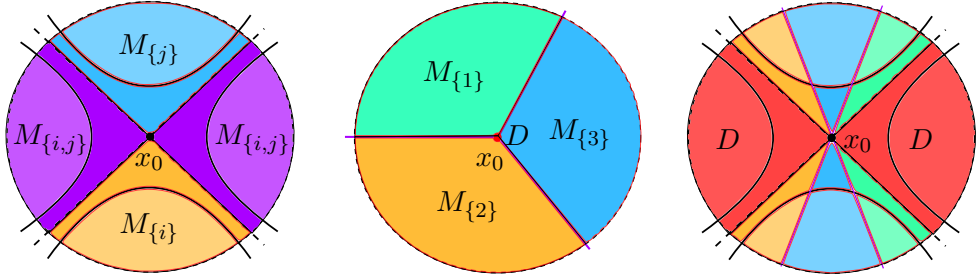

The structure of a stratified handle depends on the index of the critical point and on the stratum where it lives. There are three distinct possibilities (see Figure \ref{fig5}):

\begin{enumerate}

\item (\textit{critical point inside a regular stratum}) If $|J|=1$, the stratified handle is a smooth handle, contained in a regular stratum $M_{\{i\}}$. 

\item (\textit{critical point inside a singular stratum $M_{\{i,j\}}$}) The handle $h=j( B^{k}\times B^{m-k})$ is bisected along its core $j(B^{k}\times \{0\})$ into $h\cap M_{\{i\}}$ and $h\cap M_{\{j\}}$. By Proposition \ref{prop5}, the descending manifold of the critical point $x_0$ is contained inside $M_{\{i,j\}}$ and we obtain a bisected stratified handle, whose index is bounded by $\ind _{f}(x_0)\leq 3$.  

\item (\textit{critical point inside the singular stratum $M_{\{1,2,3\}}$}) $M_{\{1,2,3\}}$ is a smooth submanifold of codimension 2. The intersection $U_{x_0}\cap M_{\{1,2,3\}}=D$ is a 2-disk and by Proposition \ref{prop5}, the descending manifold of the critical point $x_0$ is contained in $M_{\{1,2,3\}}$. It follows that the corresponding handle is a trisected stratified handle and the index of $x_0$ is bounded above by $\ind _{f}(x_0)\leq 2$. We have three possibilities:\begin{enumerate}
\item $k=0$ and $D$ is the trisection surface of a trisected 0-handle,
\item $k=1$ and $D$ is a 2-dimensional 1-handle (whose cocore times a trisected disk gives the cocore of the 4-dimensional 1-handle),
\item $k=2$ and $D$ is the core of the 4-dimensional 2-handle, while the cocore is a trisected disk.
\end{enumerate}
\end{enumerate}

\begin{remark} In a similar fashion, we could analyze the local structure around a cri\-ti\-cal point of the function $f=\min \{f_{1},f_{2},f_{3}\}$ on $M$, which leads to the dual of handle stratification, described above. In this case, one should define a co-stratified handle as a handle whose cocore is contained inside a single stratum. Thus, we may also obtain handles with bisected core and handles with trisected core. A simple example exhibiting such handles is given by the function $f(x_1,x_2,x_3,x_4,x_5)=\min \{x_1,x_2,x_3\}$ on the 4-sphere $S^{4}=\left \{(x_{1},x_{2},x_{3},x_{4},x_{5})\in \RR ^{5}\,|\, \sum _{i=1}^{5}x_{i}^{2}=1\right \}$. 
\end{remark}

\begin{example} \label{ex2} Observe the function $f(x_{1},x_{2},x_{3},x_{4},x_{5})=\max \{x_{1},x_{2},x_{3}\}$ on the 4-sphere $$S^{4}=\left \{(x_{1},x_{2},x_{3},x_{4},x_{5})\in \RR ^{5}\,|\, \sum _{i=1}^{5}x_{i}^{2}=1\right \}\;.$$ Each regular stratum $M_{\{i\}}$ is homeomorphic to the open $4$-ball. The singular strata $M_{\{1,2\}}$, $M_{\{1,3\}}$ and $M_{\{2,3\}}$ are homeomorphic to $\inter (B^3)$ and their common boundary $M_{\{1,2,3\}}$ is homeomorphic to the 2-sphere, thus $f$ induces a $(0,0)$-trisection of $S^4$. 

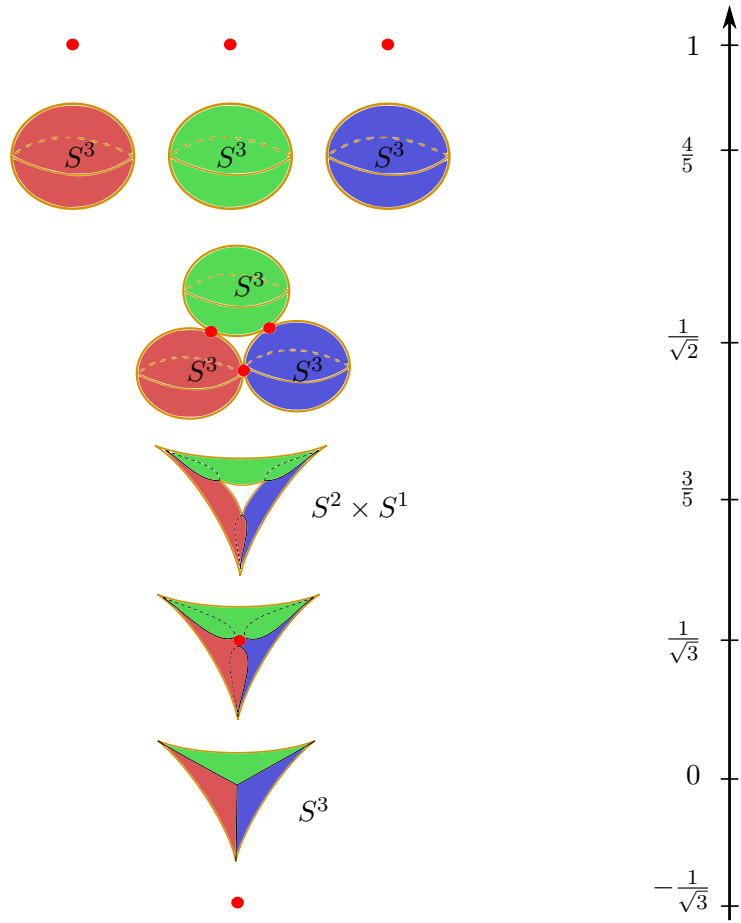
\begin{figure}[h!]
\labellist
\normalsize \hair 2pt
\pinlabel $-\frac{1}{\sqrt{3}}$ at 1150 50
\pinlabel $0$ at 1170 250
\pinlabel $S^{3}$ at 520 190
\pinlabel $\frac{1}{\sqrt{3}}$ at 1153 477
\pinlabel $\frac{1}{\sqrt{2}}$ at 1153 995
\pinlabel $\frac{3}{5}$ at 1160 740
\pinlabel $S^{2}\times S^{1}$ at 600 710
\pinlabel $\frac{4}{5}$ at 1160 1312
\pinlabel $S^3$ at 120 1312
\pinlabel $S^3$ at 380 1312
\pinlabel $S^3$ at 650 1312
\pinlabel $1$ at 1170 1500
\pinlabel $S^3$ at 330 945
\pinlabel $S^3$ at 510 945
\pinlabel $S^3$ at 410 1090
\endlabellist
\begin{center}
\includegraphics[scale=0.22]{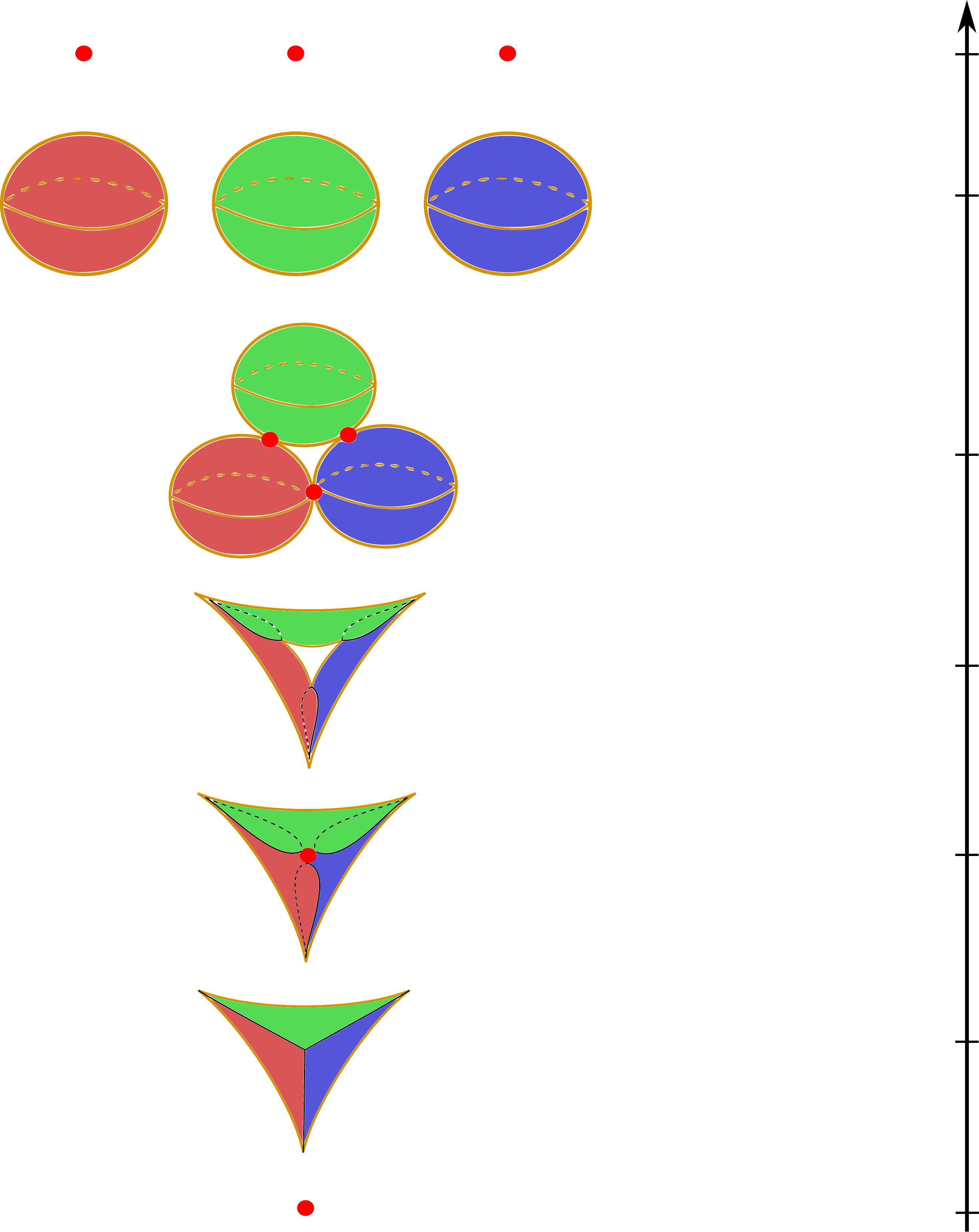}
\caption{The fibers of the CS function $f$ on $S^4$ from Example \ref{ex2}. Three different shades in every fiber indicate the decomposition into the three regular strata and critical points are represented by the red dots.}
\label{fig1}
\end{center}
\end{figure}

Let us look at the fibers of the CS Morse function $f$. Its image $f(S^{4})=\left [-\frac{1}{\sqrt{3}},1\right ]$ contains four critical values: $-\frac{1}{\sqrt{3}},\frac{1}{\sqrt{3}},\frac{1}{\sqrt{2}}$ and $1$. The critical fiber at $-\frac{1}{\sqrt{3}}$ is a single point that belongs to $M_{\{1,2,3\}}$. A regular fiber $f^{-1}(t)$ for $t\in \left (-\frac{1}{\sqrt{3}},\frac{1}{\sqrt{3}}\right )$ is homeomorphic to $S^3$:
\begin{xalignat*}{1}
& f^{-1}(0)=\left \{(0,x_{2},x_{3},x_{4},x_{5})\,|\, x_{2}^{2}+x_{3}^{2}+x_{4}^{2}+x_{5}^{2}=1, x_{2}\leq 0, x_{3}\leq 0\right \}\cup \\
& \left \{(x_{1},0,x_{3},x_{4},x_{5})\,|\, x_{1}^{2}+x_{3}^{2}+x_{4}^{2}+x_{5}^{2}=1, x_{1}\leq 0, x_{3}\leq 0\right \}\cup \\
& \left \{(x_{1},x_{2},0,x_{4},x_{5})\,|\, x_{1}^{2}+x_{2}^{2}+x_{4}^{2}+x_{5}^{2}=1, x_{1}\leq 0, x_{2}\leq 0\right \}=A_{1}\cup A_{2}\cup A_{3}
\end{xalignat*} In the above decomposition, $A_{i}$ is homeomorphic to $B^3$ and $\inter (A_{i})=f^{-1}(0)\cap M_{\{i\}}$. Moreover, $A_{i}\cap A_{j}\subseteq \cl (M_{\{i,j\}})\approx B^{2}$ and $A_{1}\cap A_{2}\cap A_{3}=\partial A_{1}\cap \partial A_{2}\cap \partial A_{3}\subseteq M_{\{1,2,3\}}$ is homeomorphic to $S^1$. The stratification induced by $f$ thus gives a 3-dimensional trisection of the regular fiber. The singular fiber $f^{-1}\left (\frac{1}{\sqrt{3}}\right )$ has a similar decomposition with a ``pinch'' in the middle (the common intersection of the three sectors is a single point instead of a circle). The pinching point $\left (\frac{1}{\sqrt{3}},\frac{1}{\sqrt{3}},\frac{1}{\sqrt{3}},0,0\right )\in M_{\{1,2,3\}}$ is a critical point of $f$. A regular fiber above a point $t\in \left (\frac{1}{\sqrt{3}},\frac{1}{\sqrt{2}}\right )$ is a union of three copies of $S^{2}\times B^{1}$, touching pairwise along the connected components of their boundaries, which is homeomorphic to $S^{2}\times S^{1}$. 

The singular fiber $f^{-1}\left (\frac{1}{\sqrt{2}}\right )$ is homeomorphic to the union of three 3-spheres, touching pairwise in three critical points $\left (\frac{1}{\sqrt{2}},\frac{1}{\sqrt{2}},0,0,0\right )$, $\left (\frac{1}{\sqrt{2}},0,\frac{1}{\sqrt{2}},0,0\right )$ and $\left (0,\frac{1}{\sqrt{2}},\frac{1}{\sqrt{2}},0,0\right )$. A regular fiber above that level is a disjoint union of three 3-spheres contained in the three regular strata, and the final critical fiber $f^{-1}(1)$ is a union of three critical points $(1,0,0,0,0)$, $(0,1,0,0,0)$ and $(0,0,1,0,0)$. See Figure \ref{fig1}. 

Each critical point of $f$ corresponds to the attachement of a stratified 4-dimensional handle, which gives a 4-manifold with boundary with a relative trisection. We start with a trisected 0-handle at the critical point $\left (-\frac{1}{\sqrt{3}},-\frac{1}{\sqrt{3}},-\frac{1}{\sqrt{3}},0,0\right )$, whose boundary is a trisected 3-sphere. The central link of this trisected $S^{3}$ is an unlink, along which a 2-handle is attached. The core of the 2-handle is a disk, lying in $M_{\{1,2,3\}}$, whose center is the second critical point $\left (\frac{1}{\sqrt{3}},\frac{1}{\sqrt{3}},\frac{1}{\sqrt{3}},0,0\right )$. The cocore of the 2-handle is a trisected disk. After the attachement of the 0-handle and the 2-handle, the boundary of the resulting manifold is homeomorphic to $S^{2}\times S^{1}$, split along three copies of $S^{2}$ by the singular strata $M_{\{i,j\}}$. At the critical fiber $f^{-1}\left (\frac{1}{\sqrt{2}}\right )$, a bisected $3$-handle is attached along the tubular neighborhood of each of these 2-spheres. After this simultaneous attachement of three 3-handles, the boundary of the resulting manifold is a disjoint union of three copies of $S^{3}$, each in their own regular stratum. At last, three smooth 4-handles are attached at the critical points $(1,0,0,0,0)$, $(0,1,0,0,0)$ and $(0,0,1,0,0)$. 

\end{example}

Our local analysis of the fibers of a CS Morse function $\max \{f_{1},f_{2},f_{3}\}$ together with Theorem \ref{th1} implies the following.

\begin{proposition} Let $f=\max \{f_{1},f_{2},f_{3}\}$ be a CS Morse function on a closed, connected and oriented 4-manifold $M$. Suppose that in each regular stratum $M_{\{i\}}$, $f$ has a single critical point of index 4, $k$ critical points of index 3 and no other critical points for $i\in \{1,2,3\}$. Moreover, suppose  that in each singular stratum $M_{\{i,j\}}$, $f$ has a single critical point of index 3, $g$ critical points of index 2 and no other critical points for all $i\neq j\in \{1,2,3\}$. Then $f$ induces a stratified handle decomposition.
\end{proposition}

Note that the stratified handles, corresponding to the critical points of $f$ inside $M_{\{1,2,3\}}$, are trisected handles. All handles, corresponding to the critical points inside other strata, come in threes: a triple of bisected handles, corresponding to critical points inside $M_{\{i,j\}}$ for $i\neq j\in \{1,2,3\}$, or a triple of smooth handles contained in $M_{\{i\}}$ for $i\in \{1,2,3\}$. A typical handle decomposition of $M$, induced by $f$, consists of the following: \\
$\left. \begin{tabular}{ll}
$\left. \begin{tabular}{lll}
a trisected 0-handle\\
$2g$ trisected 1-handles\\
a trisected 2-handle\\
\end{tabular}\right \} \Sigma _{g}\times B^{2}$\\

$3g$ bisected 2-handles\\
$3$ bisected 3-handles \\
\end{tabular}\right \} \textrm{tubular neighborhood of the spine}$\\
$\left. \begin{tabular}{ll}
$3k$ smooth 3-handles\\
$3$ smooth 4-handles \\
\end{tabular}\right \} \textrm{sector interiors}$

We have shown that the triple symmetry of trisections is not restricted to an abstract identification of the sectors, but is also a local phenomenon as we approach the trisection surface. Moreover, the symmetry may be seen on the level of handles, which might provide some new algebraic implications of trisection theory. \\ 

\begin{figure}
\labellist
\normalsize \hair 2pt
\pinlabel $-\frac{1}{\sqrt{3}}$ at 1150 50
\pinlabel $0$ at 1170 250
\pinlabel $S^{3}$ at 580 240
\pinlabel $\frac{1}{\sqrt{3}}$ at 1153 477
\pinlabel $\frac{1}{\sqrt{2}}$ at 1153 995
\pinlabel $\frac{3}{5}$ at 1160 740
\pinlabel $S^{2}\times S^{1}$ at 680 730
\pinlabel $\frac{4}{5}$ at 1160 1312
\pinlabel $S^3$ at 120 1312
\pinlabel $S^3$ at 380 1312
\pinlabel $S^3$ at 650 1312
\pinlabel $1$ at 1170 1500
\pinlabel $S^3$ at 330 945
\pinlabel $S^3$ at 510 945
\pinlabel $S^3$ at 410 1090
\endlabellist
\begin{center}
\includegraphics[scale=0.20]{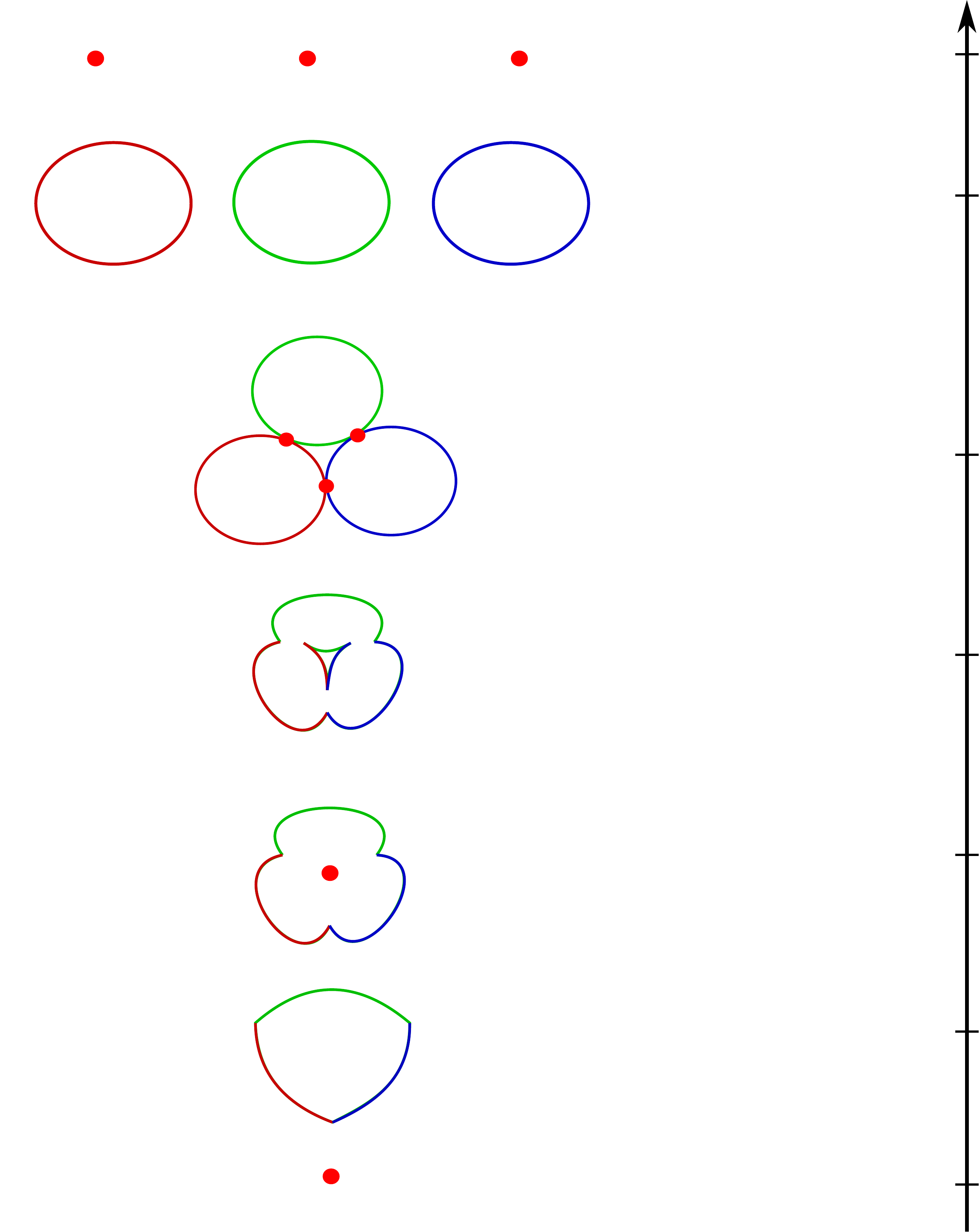}
\caption{A bridge decomposition of the unknotted 2-sphere by the CS Morse function from Example \ref{ex2}. The sphere intersects the trisection surface in two critical points $\pm \left (\frac{1}{\sqrt{3}},\frac{1}{\sqrt{3}},\frac{1}{\sqrt{3}}\right )$. The triple of critical points $\left (\frac{1}{\sqrt{2}},\frac{1}{\sqrt{2}},0\right ), \left (0,\frac{1}{\sqrt{2}},\frac{1}{\sqrt{2}}\right ), \left (\frac{1}{\sqrt{2}},0,\frac{1}{\sqrt{2}}\right )$ represent the maxima of the arcs in the triplane diagram, while the triple of critical points $(1,0,0), (0,1,0),(0,0,1)$ represent the centres of the three disks in the bridge decomposition.}
\label{fig4}
\end{center}
\end{figure}

The CS Morse function $f$ from Example \ref{ex2} may be restricted to any embedded surface in $S^{4}$. Restriction to the standardly embedded 2-sphere $S^{2}\subset S^{4}$ is also a CS Morse function and gives an analogous handle decomposition in dimension 2, see Figure \ref{fig4}. The induced trisection is called a \textbf{bridge decomposition} of the unknotted 2-sphere. Bridge decompositions of embedded surfaces were introduced by Meier and Zupan and have been extensively studied in the last few years \cite{JMMZ, MZ,MZ1}.

\end{document}